\documentclass[10pt]{article}
\usepackage{fcds,amsmath,graphicx,longtable,lscape}

\begin{document}
\title{AN ANT COLONY SYSTEM FOR TEAM ORIENTEERING PROBLEMS WITH TIME WINDOWS}

\author{Roberto MONTEMANNI$\!\!$\address{Istituto Dalle Molle di Studi sull'Intelligenza Artificiale (IDSIA), Galleria 2, CH-6928 Manno, Switzerland. Tel: +41 58 666 6667. Email: \emph{roberto@idsia.ch}.}, 
Luca Maria GAMBARDELLA$\!\!$\address{Istituto Dalle Molle di Studi sull'Intelligenza Artificiale (IDSIA), Galleria 2, CH-6928 Manno, Switzerland. Tel: +41 58 666 6663. Email: \emph{luca@idsia.ch}.}
}
\maketitle
\begin{abstract}
This paper discusses a heuristic approach for Team Orienteering Problems with Time Windows. The method we propose takes advantage of a solution model based on a hierarchic generalization of the original problem, which is combined with an Ant Colony System algorithm.

Computational results on benchmark instances previously adopted in the literature suggest that the algorithm we propose is effective in practice.
\end{abstract}


\section{Introduction} \label{intro}
Orienteering is a sport in which a competitor has to select a path from a starting point to a final destination, going through control points along the path. Each control point has a score associated with it. A travel cost is associated with each pair of control points. The competitor has to select a set of control points to be visited, so that the total score is maximized, while the total travel cost does not exceed a given threshold. As described above, the problem was first introduced as the Orienteering Problem (OP) by Tsiligrides \cite{tsi84} (the problem is sometimes referred to as the Selective Traveling Salesman Problem,  see Laporte and Martello \cite{lap90}). 

The Team Orienteering Problem (TOP) is a generalization of the OP where competitors are organized in teams. Each team component has to coordinate itself with her/his team-mates, and the objective function is to maximize the cumulative score of the team (each point can be visited by one competitor only). The TOP was firstly
studied in Butt and Cavalier \cite{but94} under the name Multiple Path Maximum Collection Problem and its current
name was introduced in Chao et al. \cite{cha96}. It is interesting to observe that, since with some elementary considerations (see Archetti et al. \cite{arc07}), it is easy to transform paths into tours, the TOP might be seen as a member of the wide family of Vehicle Routing Problems with Profits (see, for example, Feillet et al. \cite{fei95} for a survey on the Traveling Salesman Problems with Profits). For this reason, competitors are often referred to as vehicles, the team as the fleet, points as customers, and the starting/ending point as the depot.

The TOP has been recognized as a model of many
different real applications, such as the multi-vehicle version of the home fuel delivery problem (Golden et al. \cite{gol87}, \cite{gol88}), the recruiting of college football players (Butt and Cavalier \cite{but94}), the sport game of team orienteering
(Chao et al. \cite{cha96}, \cite{cha96a}), some applications of pickup or delivery services involving the use of common
carriers and private fleets (e.g. Hall and Racer \cite{hal95}) and
the service scheduling of routing technicians (Tang and Miller-Hooks \cite{tan05}).

Since the TOP is an NP-hard problem (Chao et al. \cite{cha96}), the research efforts mainly focus on heuristics and metaheuristics.  Butt and Cavalier \cite{but94} presented a greedy procedure. Chao et al. \cite{cha96a} proposed a five-step method and a heuristic
based on a stochastic algorithm. A tabu search algorithm was proposed by Tang and Miller-Hooks \cite{tan05}. Archetti et al. \cite{arc07} proposed two tabu search algorithms and a
variable neighborhood search algorithm. Ke et al. \cite{ke08} presented four methods based on Ant Colony Optimization. Not many exact algorithms have been proposed for the
TOP: a column generation based algorithm (Butt and Ryan \cite{but99}) and a branch and price based algorithm
(Boussier et al. \cite{bou07}). Several studies are devoted to the simpler OP: heuristic algorithms (see Golden et al. \cite{gol87}, \cite{gol88}, Chao et al. \cite{cha96}), exact approaches (see Fischetti et al. \cite{fis98} and Gendreau et al. \cite{gen98}) and bounding procedures (see Leifer and Rosenwein \cite{lei93} and Kataoka et al. \cite{kat98}), just to mention some of them.

Several variants of the TOP exist and one of them takes Time Window constraints into account (TOPTW). Time-window constraints are motivated by different practical situations and typically arise in those routing problems where each customer/location has to be visited within a predefined time interval specified by an earliest and a latest time into which the service has to start. In our context it is easy to see how time windows can make (some of) the applications listed before more realistic. A practical application is for companies dealing in home fuel delivery, that have to select the customers to be served in a given day (prizes are priorities in this case). The only known contributions for the Orienteering Problem with Time Windows (OPTW) are due to Kantor and Rosenwein \cite{kan92}, who described a heuristic algorithm based on an exhaustive exploration of the feasible solutions space, Mansini et al. \cite{man08}, which proposed a granular variable neighborhood search heuristic algorithm, and Righini and Salani \cite{rig08}, where some exact algorithms based on dynamic programming are discussed. Two contributions for the TOPTW, Vansteenwegen et al. \cite{van09} and Tricoire et al. \cite{tri09}, are available. They both discuss metaheuristic approaches. It is interesting to observe that some publications exist for a problem similar to the TOPTW: the Selective Vehicle Routing Problem with Time Windows (SVRPTW), where vehicles capacities and customers demands are also taken into account. Boussier et al. \cite{bou07} and Gueguen \cite{gue99} proposed exact methods for this problem.

The paper is organized as follows. Section \ref{pd} formally describes the Team Orienteering Problem with Time Windows. The Ant Colony System (ACS) algorithm we propose is described in Section \ref{acsop}. Computational results are presented and discussed in Section \ref{exp}. Conclusions are finally drawn in Section \ref{conc}.

\section{The Team Orienteering Problem with Time Windows}\label{pd}
The TOPTW is defined as follows. We are given a complete undirected graph $G = \{V, E\}$, with a positive weight $t_{ij}$ associated with each edge, representing the travel time between nodes $i$ and $j$. For each node $i \in V$ we have the following data: $p_i$ is a positive profit that is collected when the node is visited, $[a_i, b_i]$ is a time window defining the feasible arrival time at the node and $s_i$ is a non-negative service time, that is the amount of time which is spent to visit the node. Two special nodes, numbered $1$ and $n$, where $n = |V|$, are the endpoints of the path to be computed. We have $p_1 =p_n =0$, $s_1 =s_n =0$, $a_1 = a_n =0$ and $b_1 = b_n =T$, where $T$ is equal to the maximum feasible arrival time at node $n$, that is $T = \max_{i \in V \backslash \{1,n\}}\{b_i + s_i + t_{in}\}$. The TOPTW requires the computation of a set $\mathcal{P}$ of $m$ non-overlapping (apart from origin and destination) elementary paths, such that each path  $k \in \mathcal{P}$ is defined as an ordered sequence of nodes starting from node $1$ and ending at node $n$. Given a solution, we indicate with $v_i$ the arrival time at node $i$.  

We introduce a set $A$ of directed arcs such that $\forall \{i, j\} \in E$, $(i, j), (j, i) \in A$ and $t_{ij} = t_{ji}$. If we now define a binary variable $x^k_{ij} = 1 $ if arc $(i,j) \in A$ is in path $k$, $0$ otherwise, and an integer variable $z^{k}_{i}$ containing the time  customer $i$ is visited in path $k$ ($z^{k}_{i}$ does not mean anything if customer $i$ is not in path $k$, otherwise it corresponds to $v_i$), and $Q$ is an arbitrary large number, a mathematical formulation for the TOPTW is the following one:
\begin{align}
(TOPTW)	\ \ \	\  \max & \sum_{k \in \mathcal{P}} \sum_{(i,j) \in A} p_i x^k_{ij}&\label{h1} \\
\text{s.t. }    & \sum_{k \in \mathcal{P}} \sum_{j \in V} x^k_{ij} \leq 1 	& \forall i \in V \label{h2}\\
		& \sum_{j \in V} x^k_{1j} = 1 					& \forall k \in \mathcal{P} \label{h3}\\
		& \sum_{i \in V} x^k_{ih} - \sum_{j \in V} x^k_{hj} = 0 	& \forall h \in V \backslash \{1, n\},\ \forall k \in \mathcal{P} \label{h4}\\
		& \sum_{i \in V} x^k_{in} = 1 				& \forall k \in \mathcal{P} \label{h5}\\
		& z^{k}_{i} + t_{ij} + s_i - Q(1 - x^k_{ij}) \leq z^{k}_{j} 	& \forall k \in \mathcal{P} \label{h6}\\
		& a_i \leq z^{k}_{i} \leq b_i 					& \forall i \in V,\ \forall k \in \mathcal{P} \label{h7}\\
                & x^k_{ij} \in \{0, 1\} 					& \forall (i,j) \in A,\ \forall k \in \mathcal{P} \label{h8}\\
		& z^{k}_{i} \in \mathcal{Z}^+ 					& \forall i \in V,\ \forall k \in \mathcal{P} \label{h9}
\end{align}
where the objective function (\ref{h1}) maximizes the total prize collected. Constraints (\ref{h2}) impose that each customer is visited at most once; constraints (\ref{h3}), (\ref{h4}) and (\ref{h5}) impose that paths have a feasible structure; inequalities (\ref{h6}) and (\ref{h7}) deal with time windows restrictions; (\ref{h8}) and (\ref{h9}) define the domain of variables. Note that the special case where $m = |\mathcal{P}| = 1$ corresponds to the OPTW.

\section{An Ant Colony System for the \emph{TOPTW}}\label{acsop}
A brief description of the general ACS paradigm is provided in Section \ref{acs}. Section \ref{sol} is dedicated to the introduction of the Hierarchic Team Orienteering Problem with Time Windows (HTOPTW), which is functional to the algorithm we propose. The solution model adopted by our algorithm is presented in Section \ref{mod}. The remaining subsections are finally dedicated to the description of the ACS algorithm we propose to solve the HTOPTW. 
Application of an ACS algorithm
to a combinatorial optimization problem (HTOPTW in our case) requires the definition of a
constructive algorithm and possibly a local search. Accordingly, a
constructive procedure in which a set of artificial
ants builds feasible solutions to the HTOPTW is described in Section \ref{cons}.
A local search procedure specialized for the HTOPTW, that takes in input the
solutions built in the constructive phase and bring them to their local optimum, is described in Section \ref{ls}. 

\subsection{Ant Colony System} \label{acs}
The ACS algorithm is an element of
the Ant Colony Optimization (ACO) family of methods
(Dorigo et al. \cite{dor99}). These algorithms are based on a
computational paradigm inspired by real ant colonies and the way they
function. The main underlying idea was to use several constructive computational agents (simulating real ants). A dynamic memory
structure, which incorporates information on the effectiveness of
previous choices, based on the obtained results, guides the construction process of each
agent. The behavior of each single agent is therefore inspired by the
behavior of real ants.

The paradigm is based on the observation, made by ethologists, that
the medium used by ants to communicate information regarding
shortest paths to food, consists of \emph{pheromone trails}. A
moving ant lays some pheromone on the
ground, thus making a path by a trail of this substance. While an
isolated ant moves practically at random (exploration), an ant encountering a
previously laid trail can detect it and decide, with high
probability, to follow it, thus reinforcing the trail with its own
pheromone (exploitation). What emerges is a form of
\emph{auto-catalytic} process where the more the ants follow a trail,
the more attractive that trail becomes to be followed. The process
is thus characterized by a positive feedback loop, where the
probability with which an ant chooses a path increases with the
number of ants that have previously chosen the same path. The mechanism above is the inspiration for the algorithms of the \emph{ACO} family.

\subsection{The Hierarchic Team Orienteering Problem with Time Windows} \label{sol}

The input data for the HTOPTW are the same as for the TOPTW (see Section \ref{pd}).
The HTOPTW requires the computation of an ordered list of non-overlapping (apart from origin and destination) elementary paths $\mathcal{L} = (P_1, P_2, \dots, P_{|\mathcal{L}|})$, where each $P_k \in \mathcal{L}$ is defined as an ordered sequence of nodes starting from node $1$, ending at node $n$, such that $a_i \leq v_i  \leq b_i$ $\forall i \in P_k$; $\forall P_k, P_h \in \mathcal{L}$, $P_k$ and $P_h$ have in common only nodes $1$ and $n$ and $\bigcup_{P_k \in \mathcal{L}} P_k = V$. For each pair of nodes $\{i, j\}$ consecutively visited along the same $P_k \in \mathcal{L}$ we have $v_j = \max\{v_i +s_i +t_{ij}, a_j \}$, with $v_1 = 0$ on each path $P_k \in \mathcal{L}$. Finally note that we set $|\mathcal{L}|$ equal to an upper bound of the number of paths required to visit all the customers (e.g. $n$), and therefore empty tours are allowed in solutions. 

The objective function of HTOPTW can be defined as follows:
$$ \max \left( \sum_{P_k \in \mathcal{L}, k \leq m} (\sum_{(i,j) \in A} p_i x^{P_k}_{ij}) + \sum_{P_k \in \mathcal{L}, k > m} (M^{m-k} \sum_{(i,j) \in A} p_i x^{P_k}_{ij}) \right )$$
such that $a_i \leq v_i  \leq b_i$, $\forall P_k \in \mathcal{L}, \forall i \in P_k$. Note that $M$ is an arbitrary number larger than the optimal solution of the \emph{OPTW} problem defined on the same graph $G$. A complete mathematical formulation for the HTOPTW can be obtained by applying straightforward changes to the model for problem TOPTW discussed in Section \ref{pd}. 

Problem HTOPTW can be seen as a generalization of TOPTW where a set of non-overlapping tours are optimized in a hierarchic fashion. Note that if prizes are integer (as for the benchmarks we will analyze in Section \ref{etoptw}) the first $m$ tours represent the solution to the original TOPTW problem, while the remaining tours optimize over the feasible nodes that are not in the solution of the TOPTW, in a hierarchic way. 


From an optimization point of view, the rational behind the introduction of the problem HTOPTW for solving the TOPTW, is that keeping a hierarchic set of non-overlapping tours helps to have good fragments of tours placed in the $(m+1)$-st to last tours. These prepared fragments are used by the local search procedure (see Section \ref{ls}) to perform exchanges/insertions, aiming at improving the quality of the hierarchy of tours, and - as a side effect - of the first $m$ tours. The drawback of the approach we propose is that an additional overhead is necessary to optimize the $(m+1)$-st to last tours, that in fact are not part of the original objective function of the TOPTW. Computational results presented in Section \ref{exp} suggest however that it is worth paying for this overhead.

Note that a similar idea of optimizing something outside the proper solution in order to prepare good fragments to be used to improve the real solution itself, was proposed in a completely different context (and for a different optimization problem) in Lopez at al. \cite{lop98}. Again, the idea of keeping more than $m$ tours had already been considered also in Archetti et al. \cite{arc07} in the context of the TOP without time windows. In this case no hierarchic mechanism was introduced, and the best $m$ tours were basically selected as the current solution. Moreover, in \cite{arc07}, due to the different metaheuristic approaches implemented (namely, tabu search and variable neighborhood search), the two policies considered to evolve solutions (namely, swapping two customers of different tours and moving one customer from a tour to a different one) did not lead to an exploitation of the extra tours like in our case, where we move entire segments of routes at once (see Section \ref{ls}).

\subsection{Solution model} \label{mod}
The algorithm we propose uses a solution model in which each ant builds a single, giant tour (see Gambardella et al. \cite{gam99} and Irnich \cite{irn08}). A solution is represented as follows: nodes $1$ and $n$ are unified into a unique depot node, with outgoing travel times corresponding to node 1, and incoming travel times corresponding to node $n$. The depot with all its connections is then duplicated a number of times equal to the number of feasible nodes of the problem. Distances between copies of the depot are set to zero. In such a model, each time a copy of the depot is reached, tour duration is set back to zero. The use of the giant tour representation makes the HTOPTW problem closer to a  traditional traveling salesman problem.  A feasible solution is a tour that visits all the nodes exactly ones. Note that in case of too many duplicated depots, we will have dummy arcs of type \emph{(depot, depot)} visited by the solution. 

An advantage of a solution representation like the one we adopt is that the trails in direction of the duplicated depots are less attractive than in case of a single depot. This positively affects the quality of the solutions produced in the constructive procedure described in Section \ref{cons}. Moreover, as observed in Gambardella et al. \cite{gam99}, the giant tour approach makes easier to generalize and efficiently apply generic local searches procedures (see Section \ref{ls}).

\subsection{Construction phase} \label{cons}
The construction phase specializes that of the Ant Colony System algorithm (Dorigo and Gambardella
\cite{dor97}). The goal is to build feasible solutions for the HTOPTW. It generates
feasible solutions with a computational cost of order O($|V|^2$).

The procedure works as follows. Ants are sent out sequentially (not in parallel). Each ant iteratively starts from
node $1$ and adds new nodes until all nodes have been visited. When in node $i$, an ant applies a so-called transition rule, that is, it probabilistically chooses the next node $j$ from the set $F(i)$ of
feasible nodes. $F(i)$ contains all the nodes $j$ still to be
visited and feasible in terms of the time windows associated with them.

The ant in node $i$ chooses the next node $j$ to visit on the basis
of two factors: the \emph{pheromone trail} $\tau_{ij}$,
that contains a measure of how good it has been in the past to
include arc $(i, j)$ into a solution (it is the ``memory'' of the
colony) and the heuristic \emph{desirability} $\eta_{ij}$, which is defined as follows:
\begin{equation}
\eta_{ij} = \frac{p_j}{\max \left \{t_{ij}; a_j - v_i - s_i \right \} \cdot \left ( b_j - v_i - s_i  - t_{ij} \right ) + 1}\label{ff}
\end{equation}
The rationale behind equation (\ref{ff}) is that promising customers are those with a high prize, that are not far away from the current node $i$, and such that the associated time window in used in a suitable way. Note that equation (\ref{ff}) heavily relies on time windows, that therefore tend to drive desirabilities.
 
With probability $q_0$ the next node to visit is chosen as
the node $j$, $j \in F(i)$, for which the product $\tau_{ij} \cdot
\eta_{ij}$ is highest (deterministic rule), while with probability
$1-q_0$ the node $j$ is chosen with a probability given by
\[
{p_{ij}}_{j \in F(i)} = \frac{\tau_{ij} \cdot \eta_{ij}}{\sum_{l \in
F(i)} (\tau_{il} \cdot \eta_{il})}
\]

\begin{sloppypar}
(i.e., nodes connected by arcs with higher values of $\tau_{ij}~\cdot~\eta_{ij}$, $j \in F(i)$, have higher probability of being
chosen).
\end{sloppypar}

The value $q_0$ is given by $q_0 = 1 - \hat{n}/n$, where $\hat{n}$ represents the number of nodes (independently of the total number of nodes of the problem) we would like to choose using the probabilistic transition rule during each iteration of the constructive phase. Such an initialization of $q_0$ is common in ACS algorithms for Vehicle Routing-like problems. Justifications for its use can be found in Gambardella et al. \cite{gam99}.

In our implementation only the best ant, that is the ant that built the
solution with the largest profit since the beginning of the computation, is allowed to deposit a pheromone trail. The rationale is that
in this way a ``preferred route'' is memorized in the pheromone
trail matrix, and future ants will use this information to generate
new solutions in a neighborhood of this preferred route. If we refer to the most profitable solution generated since the beginning of the computation as $Tours_{Best}$, and to its profit as $Profit_{Best}$,  $\forall (i, j) \in Tours_{Best}$, we have the following formula for pheromone update:
\begin{equation}
\tau_{ij} = (1-\rho) \cdot \tau_{ij} + \rho \cdot Profit_{Best}  \label{1}
\end{equation}
where $\rho$ is a parameter regulating how strong the pheromone trace left by the best solution is.

Pheromone is also updated during solution building. In this case, it is removed from visited arcs. In other words, each ant,
when moving from node $i$ to node $j$, applies a pheromone updating
rule that causes the amount of pheromone trail on arc $(i, j)$ to
decrease. The rule is:
\begin{equation}
\tau_{ij} = (1-\psi) \cdot \tau_{ij} + \psi \cdot \tau_0 \label{2}
\end{equation}

where $\psi$ is a parameter regulating the evaporation of the pheromone trace over time, $\tau_0$ is the initial value of trails. The rationale for
using formula (\ref{2}) is that it forces ants to assure a certain variety in generated
solutions (if pheromone trail was not consumed by ants,
they would tend to generate very similar tours).

It was found that
good values for the algorithm's parameters are $\rho = \psi = 0.1$, $s =10$ and  $\tau_0 = 
(Profit_{First} \cdot n)^{-1}$, where
$Profit_{First}$ is the total prize collected by the solution with the highest profit among those generated by
the ant colony following the construction procedure described above, without using the
pheromone trails. Experience has shown these values - that are the same used in Gambardella et al. \cite{gam99} for the Vehicle Routing Problem with Time Windows - to be robust. We have also experimentally found that the algorithm is not affected by small changes in the values of these parameters. The
number of ants in the population was set to 10 (as in \cite{gam99}).

\subsection{Local search}\label{ls}
The second ingredient of the method we propose is a local search algorithm, which is run on each of the solutions produced in the construction phase and has the objective of taking them down to a local optimum. The local search procedure implemented is a specialized version of the CROSS exchange procedure described in Taillard et al. \cite{tai97}. The procedure is based on the exchange of two subchains of customers of the giant tour (see Section \ref{mod}). One of the two sub-chains can eventually be empty, implementing therefore a more traditional insertion routine. The advantage of the use of the giant tour from the local search point of view, is that the two sub-chains can be from the same tour, or from two different tours, eventually including replicas of the depot. The use of the giant tour leads then to a more general local search, that incorporates various classic local searches within a CROSS exchange framework. The maximal length of the sub-chains is, in our implementation, variable. It starts from a value defined by parameter $LS_{init}$. Every time no improvement has been found by the last  $LS_{wnd}$ generations of ants, the value of $LS_{init}$ is incremented by $LS_{step}$. The aim of this strategy is to makes the local search more and more accurate as the running time elapses. Note that too conservative values for these parameters may lead to a slow algorithm. Some preliminary tests suggested that $LS_{init} = LS_{wnd} = 3$ and $LS_{step} = 2$ are robust values, probably slightly conservative and privileging solution quality over convergence speed. We will keep these values for the experiments reported in the remainder of the paper.

Note that a kind of local search like the CROSS exchange we implemented (which does not perform inversions) is known to be very suitable for problems with time windows. In case of problems without time windows (or with very large time windows), different - and perhaps simpler - local search schemata are likely to be more suitable.

Our implementation is based on the formal theory presented in Irnich \cite{irn08}. In particular, we used a search technique where nodes for exchanges are examined in a \emph{random} order, and we move to the \emph{first} improving solution encountered. In case no improving solution is found among all the possible moves, the best non-improving solution is selected (following a strategy similar to that described in Chao et al. \cite{cha96}). This helps not being trapped into weak local optima. At most $NI$ consecutive non-improving solutions are accepted, then the local search is stopped.  It was found that $NI = 5$ is a suitable value. This setting was kept for the experiments reported in Section \ref{exp}. 

\section{Computational experiments}\label{exp}
Since in the literature the OP has often being separated from the TOP (for which, we recall, the OP is a special case), we have decided to split our computational experiments into the two classes: OPTW and TOPTW. Section \ref{eoptw} is then devoted to the OPTW, while Section \ref{etoptw} treats the TOPTW. 
 
\subsection{OPTW}\label{eoptw}
The reference benchmarks for the OPTW are those originally introduced in Righini and Salani \cite{rig08} and used in Mansini et al. \cite{man08}. We therefore adopted these benchmarks in our experiments. A brief description of them is given in Section \ref{ben}. Results are presented and discussed in Section \ref{resaaa}.

\subsubsection{Benchmarks}\label{ben}
We tested our algorithms on two classes of instances, originally adopted in Righini and Salani \cite{rig08} and Mansini et al. \cite{man08} and obtained from the well-known Solomon's data-set of VRPTW instances (problems \emph{c*\_ \!\!50}, \emph{c*\_ \!\!100}, \emph{r*\_ \!\!100} and \emph{rc*\_ \!\!100}, see Solomon \cite{sol87}) and from instances proposed by Cordeau et al. \cite{cor97} (problems \emph{pr*}) for the Multi-Depot Periodic Vehicle Routing Problem
(MDPVRP). The first data-set is composed by 29 instances obtained by considering the first 50 nodes of
Solomon's instances and by 48 instances obtained by considering the whole 100 nodes of the problems. 
The other data-set has been derived from Cordeau's
data-set (20 instances), considering all customers active in the same day. The delivery demand associated with each node in the original data-set is interpreted as the prize $p_i$ for that node. Cordeau's instances have up to 288 customers and the time windows are much wider than in Solomon's problems. 

\subsubsection{Experimental results}\label{res}\label{resaaa}
All the experiments reported in this paper have been performed on a computer equipped with a Dual AMD Opteron 250 2.4GHz processor with 4GB of RAM. The algorithms
were coded in ANSI C. A time limit of 3600 seconds has been set for each run, and statistics over five runs are reported for each instance.

Results are summarized in Table \ref{t1}. Columns have the following meaning. \emph{Problem} identifies the problems; in \emph{Known bounds} the best known bounds are reported (if a single value is present, it is the optimal solution value). Results are taken from Righini and Salani \cite{rig08} and Mansini et al. \cite{man08}. Some improved lower bounds are reported in parenthesis. They were obtained and communicated to the authors by Salani (\cite{sal08}) during the development of the present paper. The following set of columns report the results (prize collected and computation time) obtained by the Granular Variable Neighborhood Search (\emph{GVNS}) algorithm published in Mansini at al. \cite{man08}. Columns under \emph{ACS} contain the statistics (over the five runs considered) about the total prize collected in the solutions returned by the ACS algorithm and the computation times required to retrieve the best result of each run. Minimum (\emph{min}), average (\emph{avg}) and maximum (\emph{max}) prizes are reported for each problem. Note that the number of nodes visited by the best solution retrieved for each problem are also reported in brackets. Bold entries in the max prize collected column highlight new best solutions retrieved by our method, which were previously unknown (in two cases Salani \cite{sal08} matched our new best solutions, but we left our figures in bold anyway, since the results were found in parallel and independently). Entries in italics (last but one block) highlight solutions corresponding to problems for which ACS was not able to match the best known solution.

The ACS algorithm we propose was able to retrieve the best known solutions for all the Solomon's instances (first nine blocks of Table \ref{t1}). Note that the optimal solution was found in every run for most of the Solomon's problems. ACS was also able to improve 24 lower bounds over 27 instances for which a proven optimal solution is still unknown.   

Computation times indicate that our method is in general extremely fast on ``easy'' instances (for which an optimal solution is known). On the remaining instances computation times increase, but a lot of previously unknown best solutions are retrieved. This indicates that the strategy we propose, based on the solution of the HTOPTW (HOPTW in this case) pays. 

If we compare the results with those of the GVNS method (note that according to Dongarra \cite{don03b} the machine used in Mansini et al. \cite{man08} is approximately 10\% faster than the one we used), we can observe how our ACS method dominates GVNS in terms of solution quality (never worse, often better). If execution times are taken into account, we can see that ACS is also faster. This is not true for a few instances, where either the computation time is anyway comparable, or the solution found is better (and this justifies the extra time required by the algorithm).

On Cordeau's instances (last two blocks of Table \ref{t1}) the performance of the ACS algorithm is not as brilliant as on Solomon's problems. As observed in Righini and Salani \cite{rig08} (see also Section \ref{ben}), these problems are characterized by particularly wide time windows, while our algorithm - the local search component in particular - has been specialized for problems with narrow time windows (see Section \ref{ls}). Our approach is therefore a penalized on these problems. Moreover, our approach does not seem to scale up particularly well for large instances. 
This is a common drawback of algorithms based on Ant Colony System, because its construction phase tend to be too time consuming. 

We were not able to solve to optimality 3 problems over 8 instances for which the optimal solution is known. On the other hand, we were able to improve the best known lower bound for the remaining 2 instances for which a lower bound was previously known. We finally provide the first documented lower bound for the remaining 10 instances. No comparison is possible here with the GVNS algorithm since in Mansini et al. \cite{man08} no result is reported for Cordeau's instances.

A large variance in computational times can be inferred from the last two blocks of Table \ref{t1}. Once again, the presence of large time windows, coupled with the local search we implemented, might be an explanation for this.

{ \scriptsize
\begin{longtable}{|c|c|cc|ccc|ccc|}
\caption{OPTW. Computational results.} \label{t1}\\
\hline
$\!\!\!$Problem$\!\!\!$&Known&\multicolumn{2}{c|}{GVNS \cite{man08}}&\multicolumn{6}{c|}{ACS}\\
				&bounds&$\!\!\!$Prize$\!\!\!$&Sec&\multicolumn{3}{c|}{Prize}&\multicolumn{3}{c|}{Sec}\\
		       &&&&$\!\!\!$Min$\!\!\!$&Avg&$\!\!\!$$\!\!\!$Max (Nodes)$\!\!\!$$\!\!\!$&Min&Avg&Max\\ \hline
\endfirsthead
\caption[]{OPTW. Computational results (continued).}\\
\hline
$\!\!\!$Problem$\!\!\!$&Known&\multicolumn{2}{c|}{GVNS \cite{man08}}&\multicolumn{6}{c|}{ACS}\\
				&bounds&$\!\!\!$Prize$\!\!\!$&Sec&\multicolumn{3}{c|}{Prize}&\multicolumn{3}{c|}{Sec}\\
		       &&&&$\!\!\!$Min$\!\!\!$&Avg&$\!\!\!$$\!\!\!$Max (Nodes)$\!\!\!$$\!\!\!$&Min&Avg&Max\\ \hline
\endhead
\multicolumn{10}{r}{(continued on next page)}
\endfoot
\endlastfoot
$\!\!\!$ c101\_ \!\!50$\!\!\!$	&	270	&	$\!\!\!$270$\!\!\!$	&	0.19	&	$\!\!\!$270$\!\!\!$	&	270.0	&	$\!\!\!$270	(10)$\!\!\!$	&	0.06	&	0.2	&	0.31	\\
$\!\!\!$c102\_ \!\!50$\!\!\!$	&	300	&	$\!\!\!$300$\!\!\!$	&	0.21	&	$\!\!\!$300$\!\!\!$	&	300.0	&	$\!\!\!$300	(11)$\!\!\!$	&	0.09	&	0.19	&	0.26	\\
$\!\!\!$c103\_ \!\!50$\!\!\!$	&	320	&	$\!\!\!$320$\!\!\!$	&	0.75	&	$\!\!\!$320$\!\!\!$	&	320.0	&	$\!\!\!$320	(11)$\!\!\!$	&	0.59	&	0.76	&	0.9	\\
$\!\!\!$c104\_ \!\!50$\!\!\!$	&	340	&	$\!\!\!$340$\!\!\!$	&	0.77	&	$\!\!\!$340$\!\!\!$	&	340.0	&	$\!\!\!$340	(11)$\!\!\!$	&	0.02	&	0.04	&	0.05	\\
$\!\!\!$c105\_ \!\!50$\!\!\!$	&	300	&	$\!\!\!$300$\!\!\!$	&	0.25	&	$\!\!\!$300$\!\!\!$	&	300.0	&	$\!\!\!$300	(11)$\!\!\!$	&	0.07	&	0.09	&	0.13	\\
$\!\!\!$c106\_ \!\!50$\!\!\!$	&	280	&	$\!\!\!$280$\!\!\!$	&	0.27	&	$\!\!\!$280$\!\!\!$	&	280.0	&	$\!\!\!$280	(10)$\!\!\!$	&	0.01	&	0.14	&	0.35	\\
$\!\!\!$c107\_ \!\!50$\!\!\!$	&	310	&	$\!\!\!$310$\!\!\!$	&	0.15	&	$\!\!\!$310$\!\!\!$	&	310.0	&	$\!\!\!$310	(10)$\!\!\!$	&	0.01	&	0.11	&	0.19	\\
$\!\!\!$c108\_ \!\!50$\!\!\!$	&	320	&	$\!\!\!$320$\!\!\!$	&	0.28	&	$\!\!\!$320$\!\!\!$	&	320.0	&	$\!\!\!$320	(11)$\!\!\!$	&	0.33	&	0.67	&	0.79	\\
$\!\!\!$c109\_ \!\!50$\!\!\!$	&	340	&	$\!\!\!$340$\!\!\!$	&	0.24	&	$\!\!\!$340$\!\!\!$	&	340.0	&	$\!\!\!$340	(11)$\!\!\!$	&	0.39	&	0.84	&	0.97	\\ \hline
$\!\!\!$r101\_ \!\!50$\!\!\!$	&	126	&	$\!\!\!$126$\!\!\!$	&	22.44	&	$\!\!\!$126$\!\!\!$	&	126.0	&	$\!\!\!$126	(5)$\!\!\!$	&	0	&	0	&	0.01	\\
$\!\!\!$r102\_ \!\!50$\!\!\!$	&	198	&	$\!\!\!$192$\!\!\!$	&	10.34	&	$\!\!\!$198$\!\!\!$	&	198.0	&	$\!\!\!$198 	(9)$\!\!\!$	&	0.07	&	0.08	&	0.11	\\
$\!\!\!$r103\_ \!\!50$\!\!\!$	&	214	&	$\!\!\!$214$\!\!\!$	&	32.46	&	$\!\!\!$214$\!\!\!$	&	214.0	&	$\!\!\!$214	(9)$\!\!\!$	&	0.31	&	2.17	&	4.18	\\
$\!\!\!$r104\_ \!\!50$\!\!\!$	&	227	&	$\!\!\!$225$\!\!\!$	&	40.21	&	$\!\!\!$227$\!\!\!$	&	227.0	&	$\!\!\!$227	(10)$\!\!\!$	&	0.84	&	2.57	&	4.11	\\
$\!\!\!$r105\_ \!\!50$\!\!\!$	&	159	&	$\!\!\!$159$\!\!\!$	&	38.63	&	$\!\!\!$159$\!\!\!$	&	159.0	&	$\!\!\!$159	(6)$\!\!\!$	&	0.01	&	0.01	&	0.01	\\
$\!\!\!$r106\_ \!\!50$\!\!\!$	&	208	&	$\!\!\!$202$\!\!\!$	&	14.36	&	$\!\!\!$208$\!\!\!$	&	208.0	&	$\!\!\!$208	(10)$\!\!\!$	&	0.28	&	1.97	&	6.22	\\
$\!\!\!$r107\_ \!\!50$\!\!\!$	&	220	&	$\!\!\!$215$\!\!\!$	&	34.98	&	$\!\!\!$220$\!\!\!$	&	220.0	&	$\!\!\!$220	(9)$\!\!\!$	&	1.95	&	9.27	&	24.56	\\
$\!\!\!$r108\_ \!\!50$\!\!\!$	&	227	&	$\!\!\!$225$\!\!\!$	&	63.99	&	$\!\!\!$227$\!\!\!$	&	227.0	&	$\!\!\!$227	(10)$\!\!\!$	&	0.39	&	1.47	&	2.88	\\
$\!\!\!$r109\_ \!\!50$\!\!\!$	&	192	&	$\!\!\!$192$\!\!\!$	&	9.15	&	$\!\!\!$192$\!\!\!$	&	192.0	&	$\!\!\!$192	(8)$\!\!\!$	&	0.02	&	0.03	&	0.03	\\
$\!\!\!$r110\_ \!\!50$\!\!\!$	&	208	&	$\!\!\!$208$\!\!\!$	&	19.76	&	$\!\!\!$208$\!\!\!$	&	208.0	&	$\!\!\!$208	(9)$\!\!\!$	&	0	&	0.02	&	0.03	\\
$\!\!\!$r111\_ \!\!50$\!\!\!$	&	223	&	$\!\!\!$223$\!\!\!$	&	36.84	&	$\!\!\!$223$\!\!\!$	&	223.0	&	$\!\!\!$223	(9)$\!\!\!$	&	4.09	&	5.9	&	7.82	\\
$\!\!\!$r112\_ \!\!50$\!\!\!$	&	226	&	$\!\!\!$226$\!\!\!$	&	41.08	&	$\!\!\!$226$\!\!\!$	&	226.0	&	$\!\!\!$226	(10)$\!\!\!$	&	0.26	&	2.79	&	4.39	\\ \hline
$\!\!\!$rc101\_ \!\!50$\!\!\!$	&	180	&	$\!\!\!$180$\!\!\!$	&	1.19	&	$\!\!\!$180$\!\!\!$	&	180.0	&	$\!\!\!$180	(9)$\!\!\!$	&	0	&	0	&	0	\\
$\!\!\!$rc102\_ \!\!50$\!\!\!$	&	230	&	$\!\!\!$230$\!\!\!$	&	5.85	&	$\!\!\!$230$\!\!\!$	&	230.0	&	$\!\!\!$230	(10)$\!\!\!$	&	0	&	0.12	&	0.22	\\
$\!\!\!$rc103\_ \!\!50$\!\!\!$	&	240	&	$\!\!\!$240$\!\!\!$	&	9.68	&	$\!\!\!$240$\!\!\!$	&	240.0	&	$\!\!\!$240	(9)$\!\!\!$	&	0.29	&	0.7	&	0.99	\\
$\!\!\!$rc104\_ \!\!50$\!\!\!$	&	270	&	$\!\!\!$260$\!\!\!$	&	7.28	&	$\!\!\!$270$\!\!\!$	&	270.0	&	$\!\!\!$270	(10)$\!\!\!$	&	7.57	&	9.37	&	10.36	\\
$\!\!\!$rc105\_ \!\!50$\!\!\!$	&	210	&	$\!\!\!$210$\!\!\!$	&	3.29	&	$\!\!\!$210$\!\!\!$	&	210.0	&	$\!\!\!$210	(9)$\!\!\!$	&	0.19	&	0.55	&	0.93	\\
$\!\!\!$rc106\_ \!\!50$\!\!\!$	&	210	&	$\!\!\!$210$\!\!\!$	&	3.17	&	$\!\!\!$210$\!\!\!$	&	210.0	&	$\!\!\!$210	(8)$\!\!\!$	&	0	&	0	&	0	\\
$\!\!\!$rc107\_ \!\!50$\!\!\!$	&	240	&	$\!\!\!$230$\!\!\!$	&	7.35	&	$\!\!\!$240$\!\!\!$	&	240.0	&	$\!\!\!$240	(8)$\!\!\!$	&	0.02	&	3.52	& 5.05	\\
$\!\!\!$rc108\_ \!\!50$\!\!\!$	&	250	&	$\!\!\!$250$\!\!\!$	&	12.21	&	$\!\!\!$250$\!\!\!$	&	250.0	&	$\!\!\!$250	(9)$\!\!\!$	&	4.19	&	8.59	&	11.72	\\ \hline
$\!\!\!$c101\_ \!\!100$\!\!\!$	&	320	&	$\!\!\!$320$\!\!\!$	&	0.96	&	$\!\!\!$320$\!\!\!$	&	320.0	&	$\!\!\!$320		(10)$\!\!\!$	&	0.17	&	0.47	&	0.96		\\
$\!\!\!$c102\_ \!\!100$\!\!\!$	&	360	&	$\!\!\!$360$\!\!\!$	&	5.57	&	$\!\!\!$360$\!\!\!$	&	360.0	&	$\!\!\!$360		(11)$\!\!\!$	&	0.21	&	0.69	&	1.15		\\
$\!\!\!$c103\_ \!\!100$\!\!\!$	&	400	&	$\!\!\!$400$\!\!\!$	&	3.08	&	$\!\!\!$400$\!\!\!$	&	400.0	&	$\!\!\!$400		(11)$\!\!\!$	&	3.47	&	16.88	&	36.89		\\
$\!\!\!$c104\_ \!\!100$\!\!\!$	&	420	&	$\!\!\!$420$\!\!\!$	&	3.6	&	$\!\!\!$420$\!\!\!$	&	420.0	&	$\!\!\!$420		(11)$\!\!\!$	&	1.65	&	33.47	&	60.16		\\
$\!\!\!$c105\_ \!\!100$\!\!\!$	&	340	&	$\!\!\!$340$\!\!\!$	&	3.25	&	$\!\!\!$340$\!\!\!$	&	340.0	&	$\!\!\!$340		(10)$\!\!\!$	&	0.31	&	0.86	&	1.38		\\
$\!\!\!$c106\_ \!\!100$\!\!\!$	&	340	&	$\!\!\!$340$\!\!\!$	&	1.13	&	$\!\!\!$340$\!\!\!$	&	340.0	&	$\!\!\!$340		(10)$\!\!\!$	&	0.42	&	1.02	&	1.86		\\
$\!\!\!$c107\_ \!\!100$\!\!\!$	&	370	&	$\!\!\!$370$\!\!\!$	&	0.35	&	$\!\!\!$370$\!\!\!$	&	370.0	&	$\!\!\!$370		(11)$\!\!\!$	&	0.97	&	2.06	&	3.28		\\
$\!\!\!$c108\_ \!\!100$\!\!\!$	&	370	&	$\!\!\!$370$\!\!\!$	&	0.9	&	$\!\!\!$370$\!\!\!$	&	370.0	&	$\!\!\!$370		(11)$\!\!\!$	&	0.31	&	0.76	&	1.23		\\
$\!\!\!$c109\_ \!\!100$\!\!\!$	&	380	&	$\!\!\!$380$\!\!\!$	&	1.9	&	$\!\!\!$380$\!\!\!$	&	380.0	&	$\!\!\!$380		(11)$\!\!\!$	&	0.28	&	0.84	&	1.55		\\ \hline
$\!\!\!$c201\_ \!\!100$\!\!\!$	&	$\!\!\!$860(870)-1010$\!\!\!$	&	$\!\!\!$860$\!\!\!$	&	1.24	&	$\!\!\!$870$\!\!\!$	&	870.0	&	$\!\!\!$\textbf{870		(30)}$\!\!\!$	&	6.15	&	189.81	&	408.77		\\
$\!\!\!$c202\_ \!\!100$\!\!\!$	&	$\!\!\!$920-1040$\!\!\!$	&	$\!\!\!$920$\!\!\!$	&	$\!\!\!$1373.29$\!\!\!$	&	$\!\!\!$930$\!\!\!$	&	930.0	&$\!\!\!$	$\!\!\!$\textbf{930		(32)}$\!\!\!$	&	131.21	&	319.02	&	620.49		\\
$\!\!\!$c203\_ \!\!100$\!\!\!$	&	$\!\!\!$930-1000$\!\!\!$	&	$\!\!\!$930$\!\!\!$	&	$\!\!\!$1887.89$\!\!\!$	&	$\!\!\!$960$\!\!\!$	&	960.0	&	$\!\!\!$\textbf{960		(32)}$\!\!\!$	&	128.34	&	361.01	&	668.01		\\
$\!\!\!$c204\_ \!\!100$\!\!\!$	&	$\!\!\!$960-990$\!\!\!$	&	$\!\!\!$960$\!\!\!$	&	$\!\!\!$1908.02$\!\!\!$	&	$\!\!\!$970$\!\!\!$	&	970.0	&	$\!\!\!$\textbf{970		(32)	}$\!\!\!$&	415.65	&	1617.48	&	$\!\!\!$2797.30$\!\!\!$		\\
$\!\!\!$c205\_ \!\!100$\!\!\!$	&	$\!\!\!$910-1040$\!\!\!$	&	$\!\!\!$910$\!\!\!$	&	96.55	&	$\!\!\!$900$\!\!\!$	&	906.0	&	$\!\!\!$910		(31)$\!\!\!$	&	6.97	&	28.19	&	65.89	 	\\
$\!\!\!$c206\_ \!\!100$\!\!\!$	&	$\!\!\!$920-1050$\!\!\!$	&	$\!\!\!$920$\!\!\!$	&	135.51	&	$\!\!\!$920$\!\!\!$	&	920.0	&	$\!\!\!$920		(32)$\!\!\!$	&	39.29	&	87.86	&	122.00		\\
$\!\!\!$c207\_ \!\!100$\!\!\!$	&	$\!\!\!$920-1030$\!\!\!$	&	$\!\!\!$920$\!\!\!$	&	13.71	&	$\!\!\!$920$\!\!\!$	&	920.0	&	$\!\!\!$920		(32)$\!\!\!$	&	8.07	&	48.50	&	144.60		\\
$\!\!\!$c208\_ \!\!100$\!\!\!$	&	$\!\!\!$940(950)-1040$\!\!\!$ 	&	$\!\!\!$940$\!\!\!$	&	122.53	&	$\!\!\!$940$\!\!\!$	&	940.2	&	$\!\!\!$\textbf{950		(31)	}$\!\!\!$&	27.81	&	89.01	&	133.42	 	\\ \hline
$\!\!\!$r101\_ \!\!100$\!\!\!$	&	198	&	$\!\!\!$198$\!\!\!$	&	41.88	&	$\!\!\!$198$\!\!\!$	&	198.0	&	$\!\!\!$198		(9)$\!\!\!$	&	0.01	&	0.05	&	0.07		\\
$\!\!\!$r102\_ \!\!100$\!\!\!$	&	286	&	$\!\!\!$286$\!\!\!$	&	516.27	&	$\!\!\!$286$\!\!\!$	&	286.0	&	$\!\!\!$286		(11)$\!\!\!$	&	0.98	&	11.11	&	25.91		\\
$\!\!\!$r103\_ \!\!100$\!\!\!$	&	293	&	$\!\!\!$292$\!\!\!$	&	456.62	&	$\!\!\!$292$\!\!\!$	&	292.6	&	$\!\!\!$293		(11)$\!\!\!$	&	14.90	&	640.63	&	$\!\!\!$1392.33$\!\!\!$		\\
$\!\!\!$r104\_ \!\!100$\!\!\!$	&	303	&	$\!\!\!$303$\!\!\!$	&	749.09	&	$\!\!\!$303$\!\!\!$	&	303.0	&	$\!\!\!$303		(12)$\!\!\!$	&	32.47	&	164.00	&	464.73		\\
$\!\!\!$r105\_ \!\!100$\!\!\!$	&	247	&	$\!\!\!$247$\!\!\!$	&	105.13	&	$\!\!\!$247$\!\!\!$	&	247.0	&	$\!\!\!$247		(11)$\!\!\!$	&	0.27	&	3.01	&	6.38		\\
$\!\!\!$r106\_ \!\!100$\!\!\!$	&	293	&	$\!\!\!$293$\!\!\!$	&	383.55	&	$\!\!\!$293$\!\!\!$	&	293.0	&	$\!\!\!$293		(11)$\!\!\!$	&	2.82	&	86.31	&	265.91		\\
$\!\!\!$r107\_ \!\!100$\!\!\!$	&	299	&	$\!\!\!$297$\!\!\!$	&	700.88	&	$\!\!\!$293$\!\!\!$	&	294.6	&	$\!\!\!$299		(13)$\!\!\!$	&	438.23	&	922.61	&	$\!\!\!$1413.90$\!\!\!$		\\
$\!\!\!$r108\_ \!\!100$\!\!\!$	&	308	&	$\!\!\!$306$\!\!\!$	&	418.32	&	$\!\!\!$303$\!\!\!$	&	306.0	&	$\!\!\!$308		(13)$\!\!\!$	&	41.67	&	696.10	&	$\!\!\!$1666.34$\!\!\!$		\\
$\!\!\!$r109\_ \!\!100$\!\!\!$	&	277	&	$\!\!\!$277$\!\!\!$	&	87.86	&	$\!\!\!$277$\!\!\!$	&	277.0	&	$\!\!\!$277	(12)$\!\!\!$	&	6.59	&	27.98	&	51.16		\\
$\!\!\!$r110\_ \!\!100$\!\!\!$	&	284	&	$\!\!\!$283$\!\!\!$	&	289.51	&	$\!\!\!$282$\!\!\!$	&	283.2	&	$\!\!\!$284		(13)$\!\!\!$	&	27.99	&	617.57	&	$\!\!\!$1569.31$\!\!\!$		\\
$\!\!\!$r111\_ \!\!100$\!\!\!$	&	297	&	$\!\!\!$295$\!\!\!$	&	307.14	&	$\!\!\!$295$\!\!\!$	&	296.6	&	$\!\!\!$297		(12)$\!\!\!$	&	7.57	&	484.37	&	$\!\!\!$1097.63$\!\!\!$		\\
$\!\!\!$r112\_ \!\!100$\!\!\!$	&	298	&	$\!\!\!$297$\!\!\!$	&	444.52	&	$\!\!\!$295$\!\!\!$	&	297.4	&	$\!\!\!$298		(12)$\!\!\!$	&	21.77	&	947.06	&	$\!\!\!$1605.08$\!\!\!$		\\ \hline
$\!\!\!$r201\_ \!\!100$\!\!\!$	&	$\!\!\!$767(788)-916$\!\!\!$ 	&	$\!\!\!$767$\!\!\!$	&	12.58	&	$\!\!\!$795$\!\!\!$	&	795.8 &	$\!\!\!$\textbf{797		(38)}$\!\!\!$	&	355.81	&	$\!\!\!$2339.94$\!\!\!$	&	$\!\!\!$3369.30$\!\!\!$		\\
$\!\!\!$r202\_ \!\!100$\!\!\!$	&	$\!\!\!$862-1037$\!\!\!$	&	$\!\!\!$862$\!\!\!$	&	$\!\!\!$2938.39$\!\!\!$	&	$\!\!\!$895$\!\!\!$	&	899.6	&	$\!\!\!$\textbf{903		(47)	}$\!\!\!$&	389.42	&	$\!\!\!$1724.42$\!\!\!$	&	$\!\!\!$3001.66	$\!\!\!$	\\
$\!\!\!$r203\_ \!\!100$\!\!\!$	&	931-1035	&	$\!\!\!$931$\!\!\!$	&	690.92	&	$\!\!\!$983$\!\!\!$	&	989.8	&	$\!\!\!$\textbf{993		(48)	}$\!\!\!$&	1357.57	&	$\!\!\!$2641.66$\!\!\!$	&	$\!\!\!$3502.92	$\!\!\!$	\\
$\!\!\!$r204\_ \!\!100$\!\!\!$	&	1014-1138	&	$\!\!\!$1014$\!\!\!$	&	$\!\!\!$14800.28$\!\!\!$	&	$\!\!\!$1035$\!\!\!$	&	1046.4		&$\!\!\!$\textbf{1053		(53)	$\!\!\!$}&	658.87	&	$\!\!\!$1549.15$\!\!\!$	&	$\!\!\!$1964.07	$\!\!\!$	\\
$\!\!\!$r205\_ \!\!100$\!\!\!$	&	$\!\!\!$873(932)-978 $\!\!\!$	&	$\!\!\!$873$\!\!\!$	&	$\!\!\!$2837.95$\!\!\!$	&	$\!\!\!$933$\!\!\!$	&	939.2	&	$\!\!\!$\textbf{949		(46)	}$\!\!\!$&	651.50	&	1125.36	&	$\!\!\!$1740.49	$\!\!\!$	\\
$\!\!\!$r206\_ \!\!100$\!\!\!$	&	929-1033	&	$\!\!\!$929$\!\!\!$	&	$\!\!\!$4864.67$\!\!\!$	&	$\!\!\!$961$\!\!\!$	&	983.0	&	$\!\!\!$\textbf{1008		(50)	}$\!\!\!$&	$\!\!\!$1439.12$\!\!\!$	&	$\!\!\!$1671.45$\!\!\!$	&	$\!\!\!$1951.73$\!\!\!$		\\
$\!\!\!$r207\_ \!\!100$\!\!\!$	&	943-1059	&	$\!\!\!$943$\!\!\!$	&	$\!\!\!$13839.27$\!\!\!$	&	$\!\!\!$1016$\!\!\!$	&	1026.6	&	$\!\!\!$\textbf{1035		(49)}$\!\!\!$	&	645.46	&	$\!\!\!$1090.34$\!\!\!$	&	$\!\!\!$1794.28	$\!\!\!$	\\
$\!\!\!$r208\_ \!\!100$\!\!\!$	&	1029-1145	&	$\!\!\!$1029$\!\!\!$	&	$\!\!\!$5088.56$\!\!\!$	&	$\!\!\!$1048$\!\!\!$	&	1057.2	&	$\!\!\!$\textbf{1071		(55)}$\!\!\!$	&	907.35	&	1260.93	&	$\!\!\!$1723.50$\!\!\!$		\\
$\!\!\!$r209\_ \!\!100$\!\!\!$	&	879-1004	&	$\!\!\!$879$\!\!\!$	&	2453.79	&	$\!\!\!$910$\!\!\!$	&	923.4	&	$\!\!\!$\textbf{938		(45)	}$\!\!\!$&	856.86	&	1453.24	&	1878.32	$\!\!\!$	\\
$\!\!\!$r210\_ \!\!100$\!\!\!$	&	907-1003	&	$\!\!\!$907$\!\!\!$	&	6877.15	&	$\!\!\!$944$\!\!\!$	&	956.4	&	$\!\!\!$\textbf{970		(48)	}$\!\!\!$&	390.38	&	979.10	&	1666.93	$\!\!\!$	\\
$\!\!\!$r211\_ \!\!100$\!\!\!$	&	$\!\!\!$966(978)-1049$\!\!\!$ 	&	966	&	3357.44	&	$\!\!\!$1002$\!\!\!$	&	1006.8	&	$\!\!\!$\textbf{1016		(49)}$\!\!\!$	&	994.55	&	1288.08	&	1746.42	$\!\!\!$	\\ \hline
$\!\!\!$rc101\_ \!\!100$\!\!\!$	&	219	&	$\!\!\!$219$\!\!\!$	&	19.96	&	$\!\!\!$219$\!\!\!$	&	219.0	&	$\!\!\!$219		(9)	$\!\!\!$&	0.03	&	0.20	&	0.41		\\
$\!\!\!$rc102\_ \!\!100	$\!\!\!$&	266	&	$\!\!\!$266$\!\!\!$	&	38.1	&	$\!\!\!$266$\!\!\!$	&	266.0	&	$\!\!\!$266		(10)$\!\!\!$	&	4.51	&	30.91	&	43.48		\\
$\!\!\!$rc103\_ \!\!100$\!\!\!$	&	266	&	$\!\!\!$266$\!\!\!$	&	140.84	&	$\!\!\!$266$\!\!\!$	&	266.0	&	$\!\!\!$266		(10)	$\!\!\!$&	29.02	&	57.22	&	135.52		\\
$\!\!\!$rc104\_ \!\!100$\!\!\!$	&	301	&	$\!\!\!$301$\!\!\!$	&	209.2	&	$\!\!\!$301$\!\!\!$	&	301.0	&	$\!\!\!$301		(11)$\!\!\!$	&	10.22	&	29.37	&	47.26		\\
$\!\!\!$rc105\_ \!\!100$\!\!\!$	&	244	&	$\!\!\!$244$\!\!\!$	&	36.53	&	$\!\!\!$244$\!\!\!$	&	244.0	&	$\!\!\!$244		(11)$\!\!\!$	&	0.56	&	9.74	&	19.47		\\
$\!\!\!$rc106\_ \!\!100$\!\!\!$	&	252	&	$\!\!\!$250$\!\!\!$	&	38.1	&	$\!\!\!$252$\!\!\!$	&	252.0	&	$\!\!\!$252		(11)$\!\!\!$	&	5.91	&	308.64	&	520.67		\\
$\!\!\!$rc107\_ \!\!100$\!\!\!$	&	277	&	$\!\!\!$277$\!\!\!$	&	92.18	&	$\!\!\!$277$\!\!\!$	&	277.0	&	$\!\!\!$277		(10)$\!\!\!$	&	181.85	&	502.07	&	898.18		\\
$\!\!\!$rc108\_ \!\!100$\!\!\!$	&	298	&	$\!\!\!$298$\!\!\!$	&	104.33	&	$\!\!\!$298$\!\!\!$	&	298.0	&	$\!\!\!$298		(11)$\!\!\!$	&	4.27	&	207.51	&	445.34		\\ \hline
$\!\!\!$rc201\_ \!\!100$\!\!\!$	&	$\!\!\!$776(789)-1065$\!\!\!$ 	&	$\!\!\!$776$\!\!\!$	&	588.33	&	$\!\!\!$795$\!\!\!$	&	795.0	&	$\!\!\!$\textbf{795		(33)}$\!\!\!$	&	236.62	&	590.80	&	850.05		\\
$\!\!\!$rc202\_ \!\!100$\!\!\!$	&	871-1095	&	$\!\!\!$871$\!\!\!$	&	1093.06	&	$\!\!\!$931$\!\!\!$	&	931.6	&	$\!\!\!$\textbf{932		(43)}$\!\!\!$	&	770.93	&	2545.89	&	3458.20	$\!\!\!$	\\
$\!\!\!$rc203\_ \!\!100$\!\!\!$	&	$\!\!\!$907(913)-1148$\!\!\!$ 	&	$\!\!\!$907$\!\!\!$	&	1569.86	&	$\!\!\!$974$\!\!\!$	&	976.6	&	$\!\!\!$\textbf{979		(48)}$\!\!\!$	&	1175.47	&	1705.07	&	2094.52	$\!\!\!$	\\
$\!\!\!$rc204\_ \!\!100$\!\!\!$	&	$\!\!\!$1063-1282$\!\!\!$	&	$\!\!\!$1063$\!\!\!$	&	$\!\!\!$20815.24$\!\!\!$	&	$\!\!\!$1079$\!\!\!$	&	1086.2	&	$\!\!\!$\textbf{1107		(48)}$\!\!\!$	&	1185.22	&	1353.17	&	1518.42		\\
$\!\!\!$rc205\_ \!\!100$\!\!\!$	&	$\!\!\!$810(839)-1066$\!\!\!$ 	&	$\!\!\!$810$\!\!\!$	&	689.72	&	$\!\!\!$844$\!\!\!$	&	848.4	&	$\!\!\!$\textbf{855		(40)}$\!\!\!$	&	348.01	&	1549.22	&	3403.48	$\!\!\!$	\\
$\!\!\!$rc206\_ \!\!100$\!\!\!$	&	$\!\!\!$850(861)-1088$\!\!\!$ 	&	$\!\!\!$850$\!\!\!$	&	1093.93	&	$\!\!\!$881$\!\!\!$	&	884.0	&	$\!\!\!$\textbf{888		(38)}$\!\!\!$	&	305.30	&	1776.39	&	3481.45	$\!\!\!$	\\
$\!\!\!$rc207\_ \!\!100$\!\!\!$	&	$\!\!\!$879(924)-1102$\!\!\!$ 	&	$\!\!\!$879$\!\!\!$	&	2961.96	&	$\!\!\!$954$\!\!\!$	&	960.8	&	$\!\!\!$\textbf{967	(42)}	$\!\!\!$&	642.29	&	1324.46	&	1752.91	$\!\!\!$	\\
$\!\!\!$rc208\_ \!\!100	$\!\!\!$&	$\!\!\!$977(1021)-1158$\!\!\!$ 	&	$\!\!\!$977$\!\!\!$	&	8098.35	&	$\!\!\!$1002$\!\!\!$	&	1013.2	&	$\!\!\!$\textbf{1040		(45)}$\!\!\!$	&	1112.44	&	1511.37	&	1794.75	$\!\!\!$	\\ 
\hline
$\!\!\!$pr01\_ \!\!48$\!\!\!$		&	308	&-&-&	$\!\!\!$308$\!\!\!$	&	308.0	&	$\!\!\!$308	(21)$\!\!\!$	&	97.30	&	256.23	&	563.15		\\
$\!\!\!$pr02\_ \!\!96$\!\!\!$	&	404	&-&-	&	$\!\!\!$403$\!\!\!$	&	403.8	&	$\!\!\!$404		(24)$\!\!\!$	&	218.87	&	1147.78	&	2452.55	$\!\!\!$	\\
$\!\!\!$pr03\_ \!\!144$\!\!\!$		&	394	&-&-	&	$\!\!\!$394$\!\!\!$	&	394.0	&	$\!\!\!$394		(22)$\!\!\!$	&	1089.48	&	2024.26	&	3487.34	$\!\!\!$	\\
$\!\!\!$pr04\_ \!\!192$\!\!\!$&	489	&-&-	&	$\!\!\!$473$\!\!\!$	&	482.6	&	$\!\!\!$\textit{485		(24)}$\!\!\!$	&	303.29	&	1404.67	&	3422.42	$\!\!\!$	\\
$\!\!\!$pr05\_ \!\!240$\!\!\!$	&	595	&-&-	&	$\!\!\!$576$\!\!\!$	&	576.8	&	$\!\!\!$\textit{577		(31)}$\!\!\!$&	875.62	&	2075.65	&	3574.73	$\!\!\!$	\\
$\!\!\!$pr06\_ \!\!288$\!\!\!$	&	501 - ?		&-&-&	$\!\!\!$555$\!\!\!$	&	564.6	&	$\!\!\!$\textbf{567		(29)}$\!\!\!$	&	855.47	&	2199.78	&	3587.74	$\!\!\!$	\\
$\!\!\!$pr07\_ \!\!72$\!\!\!$	&	298		&-&-&	$\!\!\!$298$\!\!\!$	&	298.0	&	$\!\!\!$298		(17)	$\!\!\!$&	1.66	&	20.11	&	40.69		\\
$\!\!\!$pr08\_ \!\!144$\!\!\!$		&	463		&-&-&	$\!\!\!$461$\!\!\!$	&	462.6	&	$\!\!\!$463		(25)$\!\!\!$	&	623.64	&	2475.96	&	3553.32	$\!\!\!$	\\
$\!\!\!$pr09\_ \!\!216$\!\!\!$&	493		&-&-&	$\!\!\!$481$\!\!\!$	&	481.8	&	$\!\!\!$\textit{482		(26)}$\!\!\!$	&	1036.70	&	2318.19	&	3674.47	$\!\!\!$	\\
$\!\!\!$pr10\_ \!\!288$\!\!\!$	&	584 - ?	&-&-	&	$\!\!\!$578$\!\!\!$	&	588.4	&	$\!\!\!$\textbf{591		(30)}$\!\!\!$	&	955.54	&	2343.51	&	3042.36	$\!\!\!$	\\ \hline
$\!\!\!$pr11\_ \!\!48$\!\!\!$	&	?		&-&-&	$\!\!\!$327$\!\!\!$	&	327.8	&	$\!\!\!$\textbf{328		(20)}$\!\!\!$	&	28.42	&	1743.72	&	3036.67	$\!\!\!$	\\
$\!\!\!$pr12\_ \!\!96$\!\!\!$		&	?		&-&-&	$\!\!\!$435$\!\!\!$	&	436.4	&	$\!\!\!$\textbf{437		(24)}$\!\!\!$	&	1217.30	&	2017.56	&	3592.60	$\!\!\!$	\\
$\!\!\!$pr13\_ \!\!144$\!\!\!$	&	?		&-&-&	$\!\!\!$439$\!\!\!$	&	441.0	&	$\!\!\!$\textbf{442		(23)}$\!\!\!$	&	838.62	&	2312.73	&	3490.68	$\!\!\!$	\\
$\!\!\!$pr14\_ \!\!192$\!\!\!$		&	?		&-&-&	$\!\!\!$483$\!\!\!$	&	494.0	&$\!\!\!$\textbf{501		(28)}$\!\!\!$	&	3.37	&	23.11	&	38.31		\\
$\!\!\!$pr15\_ \!\!240$\!\!\!$	&	?		&-&-&	$\!\!\!$522$\!\!\!$	&	524.8	&	$\!\!\!$\textbf{528		(28)}$\!\!\!$	&	0.00	&	18.20	&	38.86		\\
$\!\!\!$pr16\_ \!\!288	$\!\!\!$	&	?		&-&-&	$\!\!\!$509$\!\!\!$	&	517.8	&	$\!\!\!$\textbf{525		(28)}	$\!\!\!$&	5.26	&	25.73	&	48.18		\\
$\!\!\!$pr17\_ \!\!72$\!\!\!$	&	?	&-&-	&	$\!\!\!$358$\!\!\!$	&	358.0	&	$\!\!\!$\textbf{358		(21)}$\!\!\!$	&	75.83	&	1330.85	&	2567.23	$\!\!\!$	\\
$\!\!\!$pr18\_ \!\!144$\!\!\!$		&	?		&-&-&	$\!\!\!$465$\!\!\!$	&	488.8	&	$\!\!\!$\textbf{504		(25)}	$\!\!\!$&	137.62	&	1350.04	&	3328.18	$\!\!\!$	\\
$\!\!\!$pr19\_ \!\!216$\!\!\!$	&	?		&-&-&	$\!\!\!$467$\!\!\!$	&	475.0	&	$\!\!\!$\textbf{480		(26)}$\!\!\!$	&	4.87	&	30.04	&	51.44		\\
$\!\!\!$pr20\_ \!\!288$\!\!\!$		&	?		&-&-&	$\!\!\!$547$\!\!\!$	&	552.4	&	$\!\!\!$\textbf{556		(29)}$\!\!\!$	&	5.46	&	24.57	&	42.22		\\
\hline
\end{longtable}
}
\subsection{TOPTW}\label{etoptw}

\subsubsection{Benchmarks}\label{ben2}
A set of 387 problems has been proposed in Chao et al. \cite{cha96} for the TOP (without time windows). 270 of the problems have been solved to optimality in Boussier et al. \cite{bou07}, while for the remaining 117 instances optimality has not been proven. We tested the method we propose on these benchmarks. Our method quickly converges to the optimal value for most of the problems, but for the most difficult instances (about 5\% of the whole dataset), the solutions provided by our algorithm are up to 5\% worse than those retrieved by the ACO algorithm proposed in Ke et al. \cite{ke08} (state-of-the-art method for the TOP). We recall the reader that the approach discussed in \cite{ke08}, although based on the same metaheuristic paradigm of our method, is extremely different from the algorithm we propose. The main (but not only) difference is that the idea of a hierarchic problem, useful for problems with time windows, is not considered in \cite{ke08} (where time windows are not treated). The different design of the methods  justifies the results discussed above, and shows how the design of an algorithm is strictly related to the problem under investigation. As already discussed in Section \ref{acsop}, our algorithm is optimized to exploit time windows. The hierarchic approach we adopt goes into this direction by preparing promising solution fragments that can be exchanged with fragments of the current solution characterized by similar time windows. This introduces a layer of complexity which slows down the method in case time windows are not present.

\begin{sloppypar}
Some instances exist for a problem which is slightly different from TOPTW: the Selective Vehicle Routing Problem with Time Windows, which is described in Boussier et al. \cite{bou07}. This problem has some complications with respect to TOPTW described in Section \ref{pd}. A demand is associated with each node, and vehicles have a limited capacity. Moreover, the maximum length of tours is constrained, and consequently the departure time from the depot has also to be optimized. Since our method is not designed to cope with these constraints, these instances have not been considered. In particular, dealing with the optimization of the departure times would mean to redesign both the construction phase of the ACS algorithm, and the local search component.
\end{sloppypar}

For the reasons above, we decided to carry out experiments on the same instances used for the tests summarized in Section \ref{eoptw} for the OPTW (leaving out Solomon's problems with 50 nodes). In this case we added a parameter indicating the number of vehicles available, $m$. In the test we present, we consider $2 \leq m \leq 4$. Note that also the problems originally proposed in Chao et al. \cite{cha96} for the TOP use the same range of values for parameter $m$. 

\subsubsection{Experimental results}\label{res2}
Results are summarized in Tables \ref{t2}, \ref{t3} and \ref{t4}. Columns have the same meaning as for Table \ref{t1} in Section \ref{res}, only columns \emph{Known bounds} and \emph{GVNS} are not reported here (no information is available). 

Tables \ref{t2}, \ref{t3} and \ref{t4} suggest that the algorithm we propose is suitable also for problems with more than one vehicles. No upper bound has been computed for these problems, so it is not possible to estimate the absolute quality of the solutions retrieved. However, we expect these problems to be more difficult than the  OPTW problems considered in Section \ref{res} over the same instances. The reason is that the length of each feasible solution increases as $m$ increases. 
By crosschecking the tables, it is also interesting to observe that the contribution of each additional vehicles is more and more marginal as $m$ increases. 

It is finally interesting to observe how, on some problems (typically Cordeau's instances \emph{pr*}), the time at which the best solution is retrieved approaches the time limit of 3600 seconds. This suggests that longer computation times might have allowed better results.

{  \scriptsize
\begin{longtable}{|c|ccc|ccc|}
\caption{TOPTW, $m = 2$. Computational results of the ACS algorithm.} \label{t2}\\
\hline
Problem&\multicolumn{3}{c|}{Prize}&\multicolumn{3}{c|}{Sec}\\
		       &Min&Avg&Max (Nodes)&Min&Avg&Max\\ \hline
\endfirsthead
\caption[]{TOPTW, $m = 2$. Computational results of the ACS algorithm (continued).}\\
\hline
Problem&\multicolumn{3}{c|}{Prize}&\multicolumn{3}{c|}{Sec}\\
		       &Min&Avg&Max (Nodes)&Min&Avg&Max\\ \hline
\endhead
\multicolumn{7}{r}{(continued on next page)}
\endfoot
\endlastfoot
c101\_ \!\!100	&	580	&	588.0	&	590	(21)	&	49.57	&	110.83	&	163.36	\\
c102\_ \!\!100	&	660	&	660.0	&	660	(22)	&	279.25	&	1427.36	&	3381.34	\\
c103\_ \!\!100	&	710	&	710.0	&	710	(21)	&	430.59	&	938.86	&	1777.64	\\
c104\_ \!\!100	&	750	&	754.0	&	760	(22)	&	314.14	&	1355.04	&	2048.16	\\
c105\_ \!\!100	&	640	&	640.0	&	640	(21)	&	449.59	&	1436.61	&	3038.13	\\
c106\_ \!\!100	&	620	&	620.0	&	620	(20)	&	46.75	&	104.61	&	207.26	\\
c107\_ \!\!100	&	670	&	670.0	&	670	(22)	&	24.37	&	76.02	&	140.65	\\
c108\_ \!\!100	&	680	&	680.0	&	680	(22)	&	5.80	&	95.34	&	251.49	\\
c109\_ \!\!100	&	720	&	720.0	&	720	(22)	&	282.05	&	1817.29	&	2822.51	\\ \hline
c201\_ \!\!100	&	1450	&	1452.0	&	1460	(66)	&	93.62	&	508.59	&	1551.86	\\
c202\_ \!\!100	&	1440	&	1446.0	&	1460	(66)	&	223.74	&	1568.22	&	3140.17\\
c203\_ \!\!100	&	1440	&	1446.0	&	1460	(65)	&	1846.00	&	2334.80	&	3220.08	\\
c204\_ \!\!100	&	1430	&	1436.0	&	1440	(63)	&	206.84	&	1373.22	&	2942.62	\\
c205\_ \!\!100	&	1450	&	1452.0	&	1460	(65)	&	31.12	&	925.37	&	2566.74	\\
c206\_ \!\!100	&	1460	&	1460.0	&	1460	(65)	&	678.24	&	1876.39	&	2690.34	\\
c207\_ \!\!100	&	1450	&	1452.0	&	1460	(65)	&	189.21	&	1433.98	&	3193.58	\\
c208\_ \!\!100	&	1460	&	1462.0	&	1470	(66)	&	20.38	&	1164.20	&	1605.87	\\
\hline
r101\_ \!\!100	&	349	&	349.0	&	349	(16)	&	2.06	&	36.48	&	107.15	\\
r102\_ \!\!100	&	505	&	506.6	&	508	(21)	&	797.60	&	2006.59	&	3549.10	\\
r103\_ \!\!100	&	517	&	518.8	&	520	(22)	&	459.69	&	1095.06	&	2583.94	\\
r104\_ \!\!100	&	538	&	540.0	&	544	(23)	&	1708.55	&	2291.80	&	3367.68	\\
r105\_ \!\!100	&	453	&	453.0	&	453	(20)	&	44.19	&	1037.28	&	1687.62	\\
r106\_ \!\!100	&	515	&	523.8	&	529	(21)	&	1125.75	&	2034.73	&	3320.77	\\ 
r107\_ \!\!100	&	526	&	527.4	&	529	(21) &	552.34	&	2357.58	&	3421.24	\\
r108\_ \!\!100	&	548	&	553.6	&	556	(24)	&	839.84	&	1906.24	&	2921.28	\\
r109\_ \!\!100	&	505	&	505.4	&	506	(22)	&	90.30	&	1386.61	&	2301.20	\\
r110\_ \!\!100	&	523	&	523.4	&	525	(24)	&	164.43	&	1180.21	&	2279.10	\\
r111\_ \!\!100	&	533	&	535.6	&	538	(23)	&	717.88	&	1726.17	&	3292.45	\\
r112\_ \!\!100	&	539	&	541.6	&	543	(213	&	540.79	&	1653.61	&	3122.94	\\ \hline
r201\_ \!\!100	&	1231	&	1236.0	&	1239	(70)	&	1766.97	&	2693.58	&	3533.81	\\
r202\_ \!\!100	&	1288	&	1298.6	&	1310	(78)	&	538.95	&	2482.22	&	3295.88	\\
r203\_ \!\!100	&	1341	&	1349.4	&	1358	(80)	&	2460.71	&	3057.40	&	3543.20	\\
r204\_ \!\!100	&	1397	&	1399.0	&	1404	(89)	&	1516.41	&	2713.46	&	3625.27	\\
r205\_ \!\!100	&	1336	&	1342.0	&	1346	(82)	&	1909.67	&	2593.74	&	3207.82	\\
r206\_ \!\!100	&	1367	&	1373.8	&	1381	(84)	&	2347.97	&	2954.37	&	3414.95	\\
r207\_ \!\!100	&	1367	&	1385.0	&	1400	(88)	&	1079.52	&	2624.08	&	3529.26	\\
r208\_ \!\!100	&	1416	&	1425.4	&	1433	(92)	&	1913.39	&	2828.88	&	3475.34	\\
r209\_ \!\!100	&	1344	&	1350.6	&	1361	(86)	&	2263.24	&	2988.88	&	3424.69	\\
r210\_ \!\!100	&	1351	&	1355.4	&	1360	(81)	&	1062.36	&	2423.89	&	3377.75	\\
r211\_ \!\!100	&	1396	&	1403.0	&	1411	(90)	&	2005.61	&	2726.10	&	3265.09	\\ \hline
rc101\_ \!\!100	&	427	&	427.0	&	427	(19)	&	7.32	&	27.90	&	99.56	\\
rc102\_ \!\!100	&	494	&	497.2	&	505	(20)	&	319.09	&	1496.91	&	3075.49	\\
rc103\_ \!\!100	&	506	&	510.2	&	516	(21)	&	426.05	&	1965.91	&	3533.84	\\
rc104\_ \!\!100	&	565	&	568.8	&	575	(23)		&	1583.98	&	2381.26	&	3445.15	\\
rc105\_ \!\!100	&	478	&	478.4	&	480	(21)		&	940.80	&	1676.74	&	3369.17	\\
rc106\_ \!\!100	&	480	&	480.2	&	481	(20)		&	1121.33	&	1865.07	&	3173.21	\\
rc107\_ \!\!100	&	526	&	530.2	&	534	(21)		&	231.60	&	771.48	&	2045.81	\\
rc108\_ \!\!100	&	531	&	541.6	&	550	(22)		&	194.45	&	820.99	&	2313.77	\\ \hline
rc201\_ \!\!100	&	1354	&	1365.4	&	1376	(66)	&	528.52	&	1654.27	&	3495.90	\\
rc202\_ \!\!100	&	1457	&	1464.8	&	1472	(75)	&	1076.83	&	2619.56	&	3516.05	\\
rc203\_ \!\!100	&	1523	&	1542.2	&	1573	(80)	&	882.99	&	2176.18	&	3476.60	\\
rc204\_ \!\!100	&	1608	&	1614.8	&	1622	(89)	&	1202.02	&	2561.77	&	3573.06	\\
rc205\_ \!\!100	&	1405	&	1417.4	&	1428	(73)	&	504.28	&	2115.17	&	3322.36	\\
rc206\_ \!\!100	&	1473	&	1495.6	&	1514	(78)	&	975.66	&	2105.05	&	3253.70	\\
rc207\_ \!\!100	&	1501	&	1523.0	&	1544	(80)	&	1852.18	&	2937.45	&	3424.19	\\
rc208\_ \!\!100	&	1587	&	1608.4	&	1646	(88)	&	1459.34	&	2572.30	&	3455.41	\\ \hline
pr01\_ \!\!48	&	502	&	502.0	&	502	(32)	&	28.15	&	602.21	&	1918.67	\\
pr02\_ \!\!96	& 702 &	705.2	&	714	(439)	&	289.31	&	1017.60	&	1850.88	\\
pr03\_ \!\!144	&	727	&	731.4	&	740	(40)	&	1308.75	&2083.02	&2628.40	\\
pr04\_ \!\!192	&	877	&	883.6	&	899	(50)	&	682.31	&	2439.32	&	3409.04	\\
pr05\_ \!\!240	&	1016	&	1021.8	&	1034	(55)	&	442.58	&	2278.50	&	3383.48	\\
pr06\_ \!\!288	&	977	&	987.6	&	995	(51)	&	369.49	&	1724.14	&	2648.17	\\
pr07\_ \!\!72	&	566	&	566.0	&	566	(33)	&	354.18	&	972.26	&	1669.67	\\
pr08\_ \!\!144	&	808	&	813.0	&	819	(44)	&	277.53	&	1729.99	&	2998.38	\\
pr09\_ \!\!216	&851	&	862.0	&	880	(48)	&	2628.21	&	3130.95	&	3594.25	\\
pr10\_ \!\!288	&	1053	&	1064.6	&	1078 (57)	&	2544.51	&	2918.58	&	3534.76	\\ \hline
pr11\_ \!\!48	&	544 &	544.8&	547	(37)	&60.64	&880.30&	3380.38\\
pr12\_ \!\!96	&	759&	763.4	&768	(44) &	1682.31&	2706.06	&3466.05	\\
pr13\_ \!\!144	&	788	&798.8	&816	(45)	&2700.42&	3240.72&	3564.82	\\
pr14\_ \!\!288	&	936	&	943.6	&	952	(54)	&	1240.28	&	2100.70	&	3187.72	\\
pr15\_ \!\!72	&	1090	&	1100.0	&	1120	(59)	&	2653.13	&	3188.79	&	3447.06	\\
pr16\_ \!\!144	&	1089	&	1101.0	&	1119	(58)	&	1924.75	&	2833.99	&	3562.92	\\
pr17\_ \!\!72	&	635	&646.2&	652	(39)	&64.69	&1378.00	&2951.05\\
pr18\_ \!\!144&	874	&	888.4	&	907	(47)	&	489.81	&	1960.26	&	3372.53	\\
pr19\_ \!\!216&	932	&	940.8	&	950	(54)	&	1419.52	&	2366.93	&	3152.26	\\
pr20\_ \!\!288&	1102	&	1113.0	&	1122	(62)	&	2549.27	&	3192.33	&	3372.58	\\
\hline
\end{longtable}
}
{  \scriptsize
\begin{longtable}{|c|ccc|ccc|}
\caption{TOPTW, $m = 3$. Computational results of the ACS algorithm.} \label{t3}\\
\hline
Problem&\multicolumn{3}{c|}{Prize}&\multicolumn{3}{c|}{Sec}\\
		       &Min&Avg&Max (Nodes)&Min&Avg&Max\\ \hline
\endfirsthead
\caption[]{TOPTW, $m = 3$. Computational results of the ACS algorithm (continued).}\\
\hline
Problem&\multicolumn{3}{c|}{Prize}&\multicolumn{3}{c|}{Sec}\\
		       &Min&Avg&Max (Nodes)&Min&Avg&Max\\ \hline
\endhead
\multicolumn{7}{r}{(continued on next page)}
\endfoot
\endlastfoot
c101\_ \!\!100	&	810	&	810.0	&	810	(31)	&	71.25	&	1515.58	&	3320.80	\\
c102\_ \!\!100	&	910	&	918.0	&	920	(34)	&	106.55	&	442.80	&	1121.05	\\
c103\_ \!\!100	&	960	&	972.0	&	980	(33)	&	91.88	&	1415.30	&	2230.89	\\
c104\_ \!\!100	&	990	&	1004.0	&	1020	(35)	&	476.02	&	681.70	&	856.47\\
c105\_ \!\!100	&	870	&	870.0	&	870	(33)	&	16.07	&	1010.01	&	2143.67	\\
c106\_ \!\!100	&	870	&	870.0	&	870	(32)	&	548.69	&	1132.31	&	2698.26	\\
c107\_ \!\!100	&	910	&	910.0	&	910	(34)	&	184.38	&	892.65	&	2521.55	\\
c108\_ \!\!100	&	910	&	912.0	&	920	(33)	&	82.40	&	971.24	&	3266.16	\\
c109\_ \!\!100	&	950	&	954.0	&	970	(34)	&	383.17	&	1327.54	&	2319.47	\\ \hline
c201\_ \!\!100	&	1810	&	1810.0	&	1810	(100)	&	4.43	&	259.20	&	1082.66	\\
c202\_ \!\!100	&	1780	&	1784.0	&	1790	(98)	&	297.39	&	1603.23	&	3271.71	\\
c203\_ \!\!100	&	1710	&	1738.0	&	1760	(95)	&	15.86	&	1039.94	&	3351.56	\\
c204\_ \!\!100	&	1750	&	1758.0	&	1770	(96)	&	1.12	&	1025.80	&	3526.50	\\
c205\_ \!\!100	&	1790	&	1796.0	&	1800	(99)	&	662.65	&	2007.29	&	2966.72	\\
c206\_ \!\!100	&	1780	&	1788.0	&	1800	(99)	&	235.53	&	1861.98	&	2459.64	\\
c207\_ \!\!100	&	1780	&	1784.0	&	1790	(98)	&	37.69	&	1222.95	&	3518.62	\\
c208\_ \!\!100	&	1780	&	1786.0	&	1800	(99)	&	116.61	&	1542.72	&	2783.29	\\
\hline
r101\_ \!\!100	&	481	&	481.0	&	481	(21)	&	35.79	&	240.56	&	659.12	\\
r102\_ \!\!100	&	674	&	682.0	&	691	(31)	&	1930.56	&	2435.40	&	3287.68	\\
r103\_ \!\!100	&	726	&	729.8	&	736	(33)	&	1328.00	&	2694.89	&	3592.84	\\
r104\_ \!\!100	&	750	&	760.6	&	773	(34)	&	1394.80	&	2531.35	&	3452.78	\\
r105\_ \!\!100	&	619	&	619.2	&620	(29)	&	517.98	&	795.97	&	1193.39	\\
r106\_ \!\!100	&	708	&	715.0	&	722	(31)	&	353.40	&	755.45	&	985.10	\\
r107\_ \!\!100	&	748	&	751.6	&	757	(34)	&	576.31	&	1998.89	&	2850.31	\\
r108\_ \!\!100	&	777	&	781.8	&	790	(36)	&	1343.73	&	2388.53	&	3247.86	\\
r109\_ \!\!100	&	695	&	701.8	&	710(32)	&	96.33	&	1098.83	&	2644.81	\\
r110\_ \!\!100	&	722	&	732.6	&	737	(34)	&	187.77	&	1708.47	&	3466.45	\\
r111\_ \!\!100	&	718	&	755.4	&	770	(34)	&	81.13	&	1286.72	&	3265.10	\\
r112\_ \!\!100	&	750	&	758.8	&	769	(33)	&	1346.43	&	2091.22	&	3039.20	\\ \hline
r201\_ \!\!100	&	1424	&	1428.4	&	1432	(95)	&	522.14	&	1399.97	&	2710.82	\\
r202\_ \!\!100	&	1441	&	1445.6	&	1449	(97)	&	492.80	&	1966.22	&	3451.92	\\
r203\_ \!\!100	&	1443	&	1452.4	&	1456	(99)	&	972.29	&	2564.33	&	3439.09	\\
r204\_ \!\!100	&	1458	&	1458.0	&	1458	(100)	&	5.78	&	174.66	&	385.56	\\
r205\_ \!\!100	&	1449	&	1454.8	&	1458	(100)	&	386.75	&	1830.71	&	2629.03	\\
r206\_ \!\!100	&	1455	&	1456.8	&	1458	(100)	&	15.64	&	1328.82	&	3313.74	\\
r207\_ \!\!100	&	1458	&	1458.0	&	1458	(100)	&	7.49	&	598.63	&	2466.06	\\
r208\_ \!\!100	&	1458	&	1458.0	&	1458	(100)	&	0.13	&	425.64	&	1964.77	\\
r209\_ \!\!100	&	1455	&	1456.8	&	1458	(100)	&	63.77	&	1163.74	&	2681.52	\\
r210\_ \!\!100	&	1454	&	1456.2	&	1458	(100)	&	347.09	&	1428.08	&	3400.61	\\
r211\_ \!\!100	&	1457	&	1457.8	&	1458	(100)	&	0.68	&	7.33	&	20.39	\\ \hline
rc101\_ \!\!100	&	618	&	620.4	&	621	(29)	&	16.39	&	195.67	&	389.49	\\
rc102\_ \!\!100	&	695	&	698.8	&	710 (31)	&	153.47	&	1145.53	&	2918.14	\\
rc103\_ \!\!100	&	735	&	741.4	&	747	(31)	&	1399.25	&	2753.95	&	3488.21	\\
rc104\_ \!\!100	&	793	&	813.2	&	823	(35)	&	1364.88	&	2276.97	&	3473.90	\\
rc105\_ \!\!100	&	669	&	677.8	&	682	(30)	&	132.26	&	1092.01	&	2976.23	\\
rc106\_ \!\!100	&	688	&	692.4	&	695	(30)	&	437.09	&	1748.26	&	3553.41	\\
rc107\_ \!\!100	&	742	&	748.0	&	755	(32)	&	208.06	&	1158.45	&	2288.11	\\
rc108\_ \!\!100	&	757	&	768.2	&	783	(35)	&	665.06	&	1443.67	&	2548.13	\\ \hline
rc201\_ \!\!100	&	1663	&	1672.0	&	1681	(92)	&	234.65	&	1084.91	&	2316.65	\\
rc202\_ \!\!100	&	1690	&	1701.4	&	1706	(97)	&	308.85	&	2230.74	&	3460.65	\\
rc203\_ \!\!100	&	1714	&	1719.6	&	1724	(100)	&	30.28	&	1213.17	&	3513.87	\\
rc204\_ \!\!100	&	1721	&	1722.2	&	1724	(100)	&	32.10	&	1261.26	&	3261.65	\\
rc205\_ \!\!100	&	1656	&	1672.2	&	1698	(96)	&	22.96	&	1439.79	&	2388.08	\\
rc206\_ \!\!100	&	1699	&	1712.6	&	1722	(99)	&	1089.45	&	2185.22	&	3358.51	\\
rc207\_ \!\!100	&	1702	&	1714.8	&	1722	(99)		&	1355.01	&	2553.89	&	3175.64	\\
rc208\_ \!\!100	&	1721	&	1722.2	&	1724	(100)		&	15.33	&	893.85	&	2146.41	\\ \hline
pr01\_ \!\!48	&	619	&	619.0	&	619	(43)	&	54.34	&	116.95	&	227.27	\\
pr02\_ \!\!96	&	935	&	938.4	&	942	(56)	&	1084.95	&	2333.11	&	3485.29	\\
pr03\_ \!\!144	&	984	& 991.0	&999	(60)	& 1893.28	& 876.82 &	3559.88 \\
pr04\_ \!\!192	&	1220	&	1226.4	&	1243	(69)	&	2144.70	&	2909.04	&	3561.60	\\
pr05\_ \!\!240	&	1379	&	1386.6	&	1417	(78)	&	1983.09	&	2746.66	&	3293.53	\\
pr06\_ \!\!288	&	1345	&	1359.0	&	1370	(71)	&	1538.31	&	2719.67	&	3467.55	\\
pr07\_ \!\!72	&	735	&	738.4	&	744	(49)	&	246.67	&	1330.57	&	2812.40	\\
pr08\_ \!\!144	&	1112	&	1115.0	&	1118	(59)	&	335.57	&	2544.33	&	3566.88	\\
pr09\_ \!\!216	&	1203	&	1210.0	&	1227	(69)	&	2579.76	&	3187.73	&	3467.97	\\
pr10\_ \!\!288	&	1418	&	1457.2	&	1492	(80)	&	1042.91	&	2873.12	&	3447.16	\\ \hline
pr11\_ \!\!48	&	649&	649.0	&649	(46)	&62.74	&	237.25	&	404.33	\\
pr12\_ \!\!96	&	960	&	970.6	&	985	(59)	&	2424.98	&	2883.74	&	3356.37	\\
pr13\_ \!\!144	&	1080	&	1088.8	&	1101	(68)	&	1563.33	&	2515.51	&	2803.83\\
pr14\_ \!\!192	&	1252	&	1258.6	&	1263	(74)	&	1312.17	&	2084.66	&	249.92	\\
pr15\_ \!\!240	&	1455	&	1486.8	&	1509	(83)	&	1175.35	&	2671.82	&	3404.18	\\
pr16\_ \!\!288	&	1461	&	1481.8	&	1516	(77)	&	198.84	&	2284.53	&	3410.70	\\
pr17\_ \!\!72T	&	821	&	822.6	&	832	(54)	&	564.12	&	2049.56	&	3174.96\\
pr18\_ \!\!144	&	1174	&	1209.0	&	1229	(68)		&	2199.10	&	2854.33	&	3401.15	\\ 
pr19\_ \!\!216	&	1289	&	1299.4	&	1320	(75)	&	2425.89	&	3094.80	&	3487.74	\\
pr20\_ \!\!288	&	1484	&	1491.2	&	1505	(83)	&	1739.74	&	2822.96	&	3268.58	\\
\hline
\end{longtable}
}
{  \scriptsize
\begin{longtable}{|c|ccc|ccc|}
\caption{TOPTW, $m = 4$. Computational results of the ACS algorithm.} \label{t4}\\
\hline
Problem&\multicolumn{3}{c|}{Prize}&\multicolumn{3}{c|}{Sec}\\
		       &Min&Avg&Max (Nodes)&Min&Avg&Max\\ \hline
\endfirsthead
\caption[]{TOPTW, $m = 4$. Computational results of the ACS algorithm (continued).}\\
\hline
Problem&\multicolumn{3}{c|}{Prize}&\multicolumn{3}{c|}{Sec}\\
		       &Min&Avg&Max (Nodes)&Min&Avg&Max\\ \hline
\endhead
\multicolumn{7}{r}{(continued on next page)}
\endfoot
\endlastfoot
c101\_ \!\!100	&	1010	&	1018.0	&	1020	(41)	&	105.51	&	1049.07	&	1574.07	\\
c102\_ \!\!100	&	1140	&	1142.0	&	1150	(45)	&	195.35	&	1211.26	&	2876.07	\\
c103\_ \!\!100	&	1180	&	1186.0	&	1190	(46)	&	1755.25	&	2329.56	&	2776.30	\\
c104\_ \!\!100	&	1220	&	1226.0	&	1240	(46)	&	256.95	&	1493.94	&	2746.12	\\
c105\_ \!\!100	&	1050	&	1052.0	&	1060	(43)	&	23.33	&	716.48	&	2898.51	\\
c106\_ \!\!100	&	1050	&	1058.0	&	1070	(42)	&	32.26	&	556.80	&	1188.74	\\
c107\_ \!\!100	&	1110	&	1114.0	&	1120	(44)	&	106.84	&	411.26	&	1204.96	\\
c108\_ \!\!100	&	1110	&	1112.0	&	1120	(43)	&	319.07	&	820.09	&	2339.76	\\
c109\_ \!\!100	&	1160	&	1172.0	&	1190	(46)	&	75.78	&	915.97	&	2176.46	\\ \hline
c201\_ \!\!100	&	1810	&	1810.0	&	1810	(100)	&	0.33	&	1.15	&	2.36	\\
c202\_ \!\!100	&	1810	&	1810.0	&	1810	(100)	&	0.32	&	8.87	&	36.38	\\
c203\_ \!\!100	&	1810	&	1810.0	&	1810	(100)	&	7.86	&	43.90	&	68.74	\\
c204\_ \!\!100	&	1810	&	1810.0	&	1810	(100)	&	0.47	&	1.31	&	3.10	\\
c205\_ \!\!100	&	1810	&	1810.0	&	1810	(100)	&	0.03	&	0.25	&	0.57	\\
c206\_ \!\!100	&	1810	&	1810.0	&	1810	(100)	&	0.00	&	0.13	&	0.39	\\
c207\_ \!\!100	&	1760	&	1800.0	&	1810	(100)	&	0.00	&	5.80	&	28.71	\\
c208\_ \!\!100	&	1810	&	1810.0	&	1810	(100)	&	0.02	&	0.15	&	0.48	\\
\hline
r101\_ \!\!100	&	608	&	608.0	&	608	(29)	&	9.48	&	55.13	&	203.19	\\
r102\_ \!\!100	&	817	&	825.6	&	836	(41)	&	788.47	&	1924.54	&	2676.62	\\
r103\_ \!\!100	&	893	&	902.2	&	909	(43)	&	798.91	&	2622.17	&	3532.33	\\
r104\_ \!\!100	&	931	&	944.2	&	957	(46)	&	334.64	&	2343.92	&	3192.24	\\
r105\_ \!\!100	&	761	&	766.0	&	771	(38)	&	155.49	&	838.47	&	2394.99	\\
r106\_ \!\!100	&	885	&	889.8	&	893	(41)	&	47.98	&	929.50	&	3402.78	\\ 
r107\_ \!\!100	&	930	&	932.4	&	937	(44)	&	1357.96	&	2660.90	&	3523.34	\\
r108\_ \!\!100	&	964	&	976.2	&	994	(48)	&	1038.25	&	2596.16	&	3533.89	\\
r109\_ \!\!100	&	873	&	876.6	&	879	(40)	&	157.12	&	1374.46	&	2708.17	\\
r110\_ \!\!100	&	886	&	900.0	&	908 (43)	&	370.88	&	926.27	&	1686.63	\\
r111\_ \!\!100	&	916	&	932.4	&	944	(45)	&	855.90	&	1896.35	&	2866.03	\\
r112\_ \!\!100	&	940	&	947.6	&	954	(45)	&	900.41	&	1662.61	&	2496.41	\\ \hline
r201\_ \!\!100	&	1458	&	1458.0	&	1458	(100)	&	127.45	&	376.41	&	827.34	\\
r202\_ \!\!100	&	1458	&	1458.0	&	1458	(100)	&	182.35	&	936.24	&	2728.65	\\
r203\_ \!\!100	&	1458	&	1458.0	&	1458	(100)	&	0.02	&	73.92	&	352.55	\\
r204\_ \!\!100	&	1458	&	1458.0	&	1458	(100)	&	0.00	&	0.01	&	0.02	\\
r205\_ \!\!100	&	1458	&	1458.0	&	1458	(100)	&	0.00	&	1.31	&	5.59	\\
r206\_ \!\!100	&	1458	&	1458.0	&	1458	(100)	&	0.00	&	0.01	&	0.06	\\
r207\_ \!\!100	&	1458	&	1458.0	&	1458	(100)	&	0.00	&	0.00	&	0.00	\\
r208\_ \!\!100	&	3456	&	3456.0	&	1458	(100)	&	0.00	&	0.01	&	0.03	\\
r209\_ \!\!100	&	1458	&	1458.0	&	1458	(100)	&	0.00	&	0.00	&	0.00	\\
r210\_ \!\!100	&	1458	&	1458.0	&	1458	(100)	&	0.05	&	3.13	&	11.16	\\
r211\_ \!\!100	&	1458	&	1458.0	&	1458	(100)	&	0.00	&	0.00	&	0.00	\\ \hline
rc101\_ \!\!100	&	801	&	805.4	&	808	(37)	&	84.48	&	1324.33	&	3199.64	\\
rc102\_ \!\!100	&	896	&	899.4	&	903	(42)	&	1210.62	&	2218.81	&	3479.62	\\
rc103\_ \!\!100	&	937	&	941.8	&	948	(42)	&	288.47	&	2005.21	&	3160.23	\\
rc104\_ \!\!100	&	982	&	1013.2	&	1052	(47)	&	1331.32	&	2139.30	&	3140.49	\\
rc105\_ \!\!100	&	860	&	867.2	&	875	(40)	&	134.31	&	1052.30	&	2201.45	\\
rc106\_ \!\!100	&	888	&	901.4	&	908	(39)	&	559.05	&	2106.67	&	2892.46	\\
rc107\_ \!\!100	&	952	&	959.2	&	964	(45)	&	104.30	&	1763.11	&	3285.91	\\
rc108\_ \!\!100	&	995	&	1000.2	&	1007	(47)	&	746.11	&	2222.29	&	3534.34	\\ \hline
rc201\_ \!\!100	&	1724	&	1724.0	&	1724	(100)	&	854.00	&	2327.40	&	3427.79	\\
rc202\_ \!\!100	&	1724	&	1724.0	&	1724	(100)	&	27.23	&	1095.77	&	2801.97	\\
rc203\_ \!\!100	&	1724	&	1724.0	&	1724	(100)	&	0.00	&	308.58	&	1463.78	\\
rc204\_ \!\!100	&	1724	&	1724.0	&	1724	(100)	&	0.00	&	53.33	&	266.67	\\
rc205\_ \!\!100	&	1724	&	1724.0	&	1724	(100)	&	360.01	&	1313.73	&	2126.78	\\
rc206\_ \!\!100	&	1722	&	1723.6  &	1724	(100)	&	0.13	&	2.79	&	6.45	\\
rc207\_ \!\!100	&	1722	&	1722.4	&	1724	(100)	&	1.28	&	72.12	&	172.67	\\
rc208\_ \!\!100	&	1724	&	1724.0	&	1724	(100)	&	0.00	&	0.00	&	0.00	\\ \hline
pr01\_ \!\!48	&	654	&	655.2	&	657	(48)	&	0.52	&	15.43	&	50.61	\\
pr02\_ \!\!96	&	1063	&	1067.8	&	1072	(68)	&	1185.62	&	2164.96	&	3146.01	\\
pr03\_ \!\!144	&	1190	&	1204.6	&	1222	(77)	&	817.89	&	2073.95	&	3476.42	\\
pr04\_ \!\!192	&	1503	&	1506.0	&	1515	(85)	&	2255.22	&	2986.02	&	3414.45	\\
pr05\_ \!\!240	&	1710	&1722.2	&1740	(94)	&1710.87	&2818.47	&3528.73\\
pr06\_ \!\!288	&	1709	&	17235.4	&	1740	(95)	&	1673.42	&	2935.83	&	3445.90	\\
pr07\_ \!\!72	&	860	&	865.4	&	872	(58)	&	35.40	&	1940.19	&	3318.06	\\
pr08\_ \!\!144	&	1337	&	1351.6	&	1376	(79)	&	1247.71	&	2804.25	&	3578.14	\\
pr09\_ \!\!216	&	1544	&	1554.2	&	1561	(91)	&	3174.85	&	3371.92	&	3586.72	\\
pr10\_ \!\!288	&	1813	&	1819.8	&	1827	(101)	&	2415.82	&	3365.93	&	3456.00	\\ \hline
pr11\_ \!\!48	&	657	&	657.0	&	657	(48)	&	0.00	&	0.09	&	0.35	\\
pr12\_ \!\!96	&	1110	&	1114.0	&	1118	(74)	&	1736.62	&	2486.72	&	3002.43	\\
pr13\_ \!\!144	&1303	&	1316.4	&	1329	(85)	&	2039.13	&	2926.28	&	3399.21\\
pr14\_ \!\!192 & 	1521	&	1534.4	&	1568	(91)	&	1810.47	&	2730.20	&	3268.35	\\
pr15\_ \!\!240&		1800	&	1822.8	&	1854	(102)	&	1820.47	&	2631.27	&	3445.79	\\
pr16\_ \!\!288&		1833	&	1867.2	&	1887	(104)	&	1140.65	&	2614.66	&	3477.03	\\
pr17\_ \!\!72&		923	&	923.8	&	925	(65)	&	1557.48	&	2749.67	&	3450.60	\\
pr18\_ \!\!144&		1455	&	1464.8	&	1470	(82)	&	2621.62	&	3209.13	&	3478.46	\\
pr19\_ \!\!216&		1567	&	1579.4	&	1596	(96)	&	2740.67	&	3293.94	&	3586.13	\\
pr20\_ \!\!288&		1811	&	1821.4	&	1841	(98)	&	2661.58	&	3193.02	&	3472.65
	\\ \hline
\end{longtable}
}
\section{Conclusions}\label{conc}
An algorithm based on the Ant Colony System paradigm has been discussed for the Team Orienteering Problem with Time Windows. The method is based on the solution of a hierarchic generalization of the original Team Orienteering Problem.

Experimental results on benchmarks instances previously adopted in the literature confirm the practical effectiveness of the algorithm we propose. 
\\ $ $ \\$ $ \\
\textbf{\large Acknowledgement}\\ $ $ \\ 
The authors are indebted to Matteo Salani and Giovanni Righini for the useful discussions and for having provided the benchmark instances used in the experiments.

\bibliographystyle{plain}
\bibliography{op}
\end{document}